\documentclass[a4paper,10pt,reqno]{amsart}
\usepackage{verbatim}
\usepackage{amssymb}

\usepackage{enumerate}
\usepackage[active]{srcltx}
\numberwithin{equation}{section}

\usepackage{t1enc}
\usepackage[utf8x]{inputenc}

\newtheorem{theorem}{Theorem}[section]
\newtheorem{lemma}[theorem]{Lemma}
\newtheorem{problem}[theorem]{Problem}
\newtheorem{example}[theorem]{Example}

\newtheorem{corollary}[theorem]{Corollary}

\theoremstyle{definition}
\newtheorem{definition}[theorem]{Definition}

\theoremstyle{remark}

\newcommand{\dom}{dominating pair}
\newcommand{\pd}{pinning down pair}

\newcommand{\setm}{\setminus}

\newcommand{\subs}{\subset}
\newcommand{\ran}{\operatorname{ran}}

\newcommand{\cof}{\operatorname{cof}}

\def\<{\left\langle}
\def\>{\right\rangle}
\def\br#1;#2;{\bigl[ {#1} \bigr]^ {#2} }

\author[I. Juh\'asz]{Istv\'an Juh\'asz}
\address      { Alfr\'ed Rényi Institute of Mathematics}
\email{juhasz@renyi.hu}

\author[L. Soukup]{Lajos Soukup}
\thanks
  {
   }
\address
      { Alfr{\'e}d R{\'e}nyi Institute of Mathematics%
}
\email{soukup@renyi.hu}

\author[Z. Szentmikl\'ossy]{Zolt\'an Szentmikl\'ossy}
\address{E\"otv\"os University of Budapest}
\email{szentmiklossyz@gmail.com}

\subjclass[2010]{54A25, 54A35, 54A65, 03E35}
\keywords{density of a topological space, cardinal function, dominating pair, pinning down pair}

\title[Dominating and pinning down pairs]{Dominating and pinning down pairs for topological spaces}
\thanks{The research on and preparation of this paper was
supported by  NKFIH grant no.  K 129211.}
\date{\today}

\begin{document}

\hfill \emph{This paper is dedicated to the memory of Phil Zenor}

\bigskip

\begin{abstract}
We call a pair of infinite cardinals $(\kappa,\lambda)$ with $\kappa > \lambda$ a {\em dominating} (resp. {\em pinning down}) pair
for a topological space $X$ if for every subset $A$ of $X$ (resp. family $\mathcal{U}$ of non-empty open sets in $X$) of cardinality
$\le \kappa$ there is $B \subs X$  of cardinality $\le \lambda$ such that $A \subs \overline{B}$ (resp. $B \cap U \ne \emptyset$
for each $U \in \mathcal{U}$). Clearly, a dominating pair is also a pinning down pair for $X$.
Our definitions generalize the concepts introduced in \cite{GTW} resp. \cite{BT} which focused on pairs of the form $(2^\lambda,\lambda)$.

The main aim of this paper is to answer a large number of the numerous problems from  \cite{GTW} and \cite{BT}
that asked if certain conditions on a space $X$ together with the assumption that $(2^\lambda,\lambda)$ or $((2^\lambda)^+,\lambda)$
is a \pd\ or \dom\ for $X$ would imply $d(X) \le \lambda$.
\end{abstract}

\maketitle

\section{Introduction}

\begin{definition}\label{df:parok}
Let $X$ be a topological space and $\kappa > \lambda \ge \omega$ be cardinals.
\begin{enumerate}[(i)]
\item $(\kappa,\lambda)$ is a \pd\ for $X$ if for every
family $\mathcal{U}$ of non-empty open sets in $X$ with
$|\mathcal{U}| \le \kappa$ there is $B \subs X$  of cardinality $\le \lambda$ such that
$B \cap U \ne \emptyset$ for each $U \in \mathcal{U}$. (I.e., $\mathcal{U}$ may be pinned down with $\le \lambda$ points.)

\smallskip

\item $(\kappa,\lambda)$ is a \dom\ for $X$ if for every $A \in [X]^{\le \kappa}$
there is $B \in [X]^{\le \lambda}$ such that $A \subs \overline{B}$.
\end{enumerate}
\end{definition}

These definitions generalize the concepts introduced in \cite{GTW} and \cite{BT}, where
the focus was on pairs of the form $(2^\lambda,\lambda)$, hence the terminology was different.

Clearly, a dominating pair is also a pinning down pair for $X$. That the latter is much weaker than
the former is well illustrated by the facts that while $(\mathfrak{c},\omega)$ is a \pd\ for
every Cantor cube, see \cite{BT}, if $(\mathfrak{c},\omega)$ is a \dom\ for a compact $T_2$-space $X$
then $X$ is separable by \cite{GTW}.

If $\kappa >  d(X)$
then trivially $(\kappa, d(X))$ is a \dom\ for $X$. This leads to the following natural question:
Under what conditions does the fact that $(\kappa,\lambda)$ is a \dom\ (or \pd)\ for $X$
imply $d(X) \le \lambda$?  This, of course, was the central topic of the papers \cite{GTW} and \cite{BT}
as well, restricted to pairs of the form $(2^\lambda,\lambda)$. We hope to convince the reader
about the usefulness of our more general approach
by strengthening a number of results from and
solving quite a few of the open problems raised in \cite{GTW} and \cite{BT}.

The terminology and notation we use in this note is standard, as e.g. in \cite{Ju}.
For any topological space $X$ we shall use $\tau(X)$ to denote its topology, i.e. the family of all
open sets in $X$, moreover $\tau^+(X) = \tau(X) \setm \{\emptyset\}$. Similarly, $RO(X)$ (resp. $RC(X)$) denotes
the family of all regular open (resp. regular closed) sets in $X$ and $RO^+(X) = RO(X) \setm \{\emptyset\}$
(resp. $RC^+(X) = RC(X) \setm \{\emptyset\}$).

We shall use $\mathbb{C}_\kappa$ to denote the Cantor cube $\{0,1\}^\kappa$. This is to distinguish it
from the cardinal number $2^\kappa$.

We emphasize that,
unlike in \cite{GTW} and \cite{BT}, unless stated explicitly, no separation axiom is assumed about
the spaces occurring in our results.

\section{Results for \pd s}

\begin{theorem}\label{tm:card}
Assume that $(\kappa,\lambda)$ is a \pd\ for $X$, moreover $$|X| \le \kappa = \cof \big([\kappa]^{\lambda}, \subs \big).$$
Then $d(X) \le \lambda$.
\end{theorem}

\begin{proof}
Let us start by fixing an elementary submodel $M$ of
$H(\vartheta)$ for a large enough regular cardinal $\vartheta$  with $|M| = \kappa$ and $\{X\} \cup (\kappa + 1) \subs M$.
Note that then $|X| \le \kappa = |M|$ implies $X \subs M$.

By elementarity
and our assumption there is a $\kappa$-sized cofinal subset $\mathcal{A}$ of the partial order
$\big([\kappa]^{\lambda}, \subs \big)$ with $\mathcal{A} \in M$ and hence $\mathcal{A} \subs M$ as well.

Since $|\tau^+(X) \cap M| \le \kappa$ and $(\kappa,\lambda)$ is a \pd\ for $X$, moreover $\mathcal{A}$ is
cofinal in $\big([\kappa]^{\lambda}, \subs \big)$, there is a set $A \in \mathcal{A} \subs M$ that meets
every member of $\tau^+(X) \cap M$. But then by elementarity we have $A \cap U \ne \emptyset$ for all
$U \in \tau^+(X)$ as well, i.e. $A$ is dense in $X$.
\end{proof}

The following immediate corollary of Theorem \ref{tm:card} gives an affirmative answer to problem 4.9 of \cite{BT}.

\begin{corollary}\label{cr:4.9.BT}
If $|X| \le \mathfrak{c}$ and $(\mathfrak{c},\omega)$ is a \pd\ for $X$ then $X$ is separable.
\end{corollary}

Since hereditarily Lindelöf $T_2$-spaces have cardinality $\le \mathfrak{c}$, this also
yields an affirmative answer to problem 4.10 of \cite{BT}.

Problem 4.8 of \cite{BT} asks if the statement of Corollary \ref{cr:4.9.BT} remains valid if the condition $|X| \le \mathfrak{c}$ in it
is weakened to $d(X) \le \mathfrak{c}$. This however is simply false. Indeed, as we have noted above,
$(\mathfrak{c},\omega)$ is a \pd\ for
the Cantor cube $\mathbb{C}_{2^\mathfrak{c}}$ of weight ${2^\mathfrak{c}}$,  moreover $\mathbb{C}_{2^\mathfrak{c}}$
is known to have density $\log{2^\mathfrak{c}}$,
but clearly $\omega < \log{2^\mathfrak{c}} \le \mathfrak{c}$.

It is also well-known that if $\lambda < \kappa < \aleph_\omega$ then $ \kappa = \cof \big([\kappa]^{\lambda}, \subs \big)$ holds.
This then yields the following corollary of  Theorem \ref{tm:card}.

\begin{corollary}\label{cr:alefw}
If $\,|X| < \aleph_\omega$ and $(|X|,\lambda)$ is a \pd\ for the space $X$ then $d(X) \le \lambda$.
\end{corollary}

The {\em pinning down number} $pd(X)$ of a space $X$ is defined as the smallest cardinal $\lambda$
such that for every neighborhood assignment $\{U_x : x \in X\}$ there is a set $B$ of size $\le \lambda$
with $B \cap U_x \ne \emptyset$ for all $x \in X$. It is clear that if $\Delta(X) = |X| > pd(X)$ then
$(|X|,pd(X))$ is a \pd\ for the space $X$. (Here $\Delta(X) = \min \{|U| : U \in \tau^+(X)\}$.)
So it follows from Corollary \ref{cr:alefw} that for any space $X$ with $\Delta(X) = |X| < \aleph_\omega$
we have $d(X) = pd(X)$. This in turn, by Lemma 2.2 of \cite{JSSz}, implies $d(X) = pd(X)$ whenever $|X| < \aleph_\omega$,
a result first proved in Theorem 5.1 of \cite{BR}.

\medskip

Before formulating our next result we recall that a space $X$ is called {\em quasiregular} if every non-empty open
set in $X$ includes a non-empty regular closed set. Moreover, $F(X)$ denotes the supremum of the sizes of free subsets
of $X$, where a set is free if it has a well-ordering that turns it into a free sequence, see \cite{JSSz1}.

\begin{theorem}\label{tm:qr}
Assume that $X$ is a topological space and $\lambda$ is a cardinal such that
for every set $S \in [X]^\lambda$ there is $R \in RC^+(X)$ with $R \cap \overline{S} = \emptyset$,
moreover $(2^\lambda,\lambda)$ is a \pd\ for $X$.
Then  $F(X) > \lambda$, i.e. there is a free sequence of length $\lambda^+$ in $X$.
\end{theorem}

\begin{proof}
Let $M$ be a $\lambda$-closed elementary submodel of
$H(\vartheta)$ for a large enough regular cardinal $\vartheta$  with $|M| = 2^\lambda$ and $X \in M$.
Let us put $\mathcal{R} = M \cap RC^+(X)$, then $|\mathcal{R}| \le 2^\lambda$ implies the existence of
a set $P \in [X]^\lambda$ that pins down $\mathcal{R}$, i.e. $P \cap R \ne \emptyset$ for all $R \in \mathcal{R}$.
Also, by elementarity, for each set $S \in [X \cap M]^\lambda$  there is $R \in  \mathcal{R}$ with
$R \cap \overline{S} = \emptyset$.

We claim that there is a point $p \in P$ such that for each $S \in [X \cap M]^\lambda$  there is $R \in \mathcal{R}_p$ with
$R \cap \overline{S} = \emptyset$, where $\mathcal{R}_p = \{R \in \mathcal{R} : p \in R\}$.
Indeed, assume otherwise, then for every point $p \in P$ there is a set $S_p \in [X \cap M]^\lambda$ such that
$R \cap \overline{S_p} \ne \emptyset$ for all $R \in \mathcal{R}_p$. But then for the set $S = \bigcup_{p \in P} S_p \in [X \cap M]^\lambda$
we would have $R \cap \overline{S} \ne \emptyset$ for all $R \in \mathcal{R}$, which is a contradiction.

So, let us fix $p \in P$ as in the claim and put $\mathcal{W} = \{X \setm R : R \in \mathcal{R}_p \}$.
Then the set $A = X \cap M$, the family $\mathcal{W}$ of open sets and the cardinal $\lambda^+$
satisfy the requirements of Lemma 2.1 of \cite{JSSz1}, namely (a) $\mathcal{W}$
is closed under unions of subfamilies of size $\lambda$, (b) $A \setm W \ne \emptyset$ for all $W \in \mathcal{W}$,
and (c) for any $S \in [A]^\lambda$ there is $W \in \mathcal{W}$ with $\overline{S} \subs W$.
Consequently, by the lemma there is a free sequence of length $\lambda^+$ in $X$ included in $A$.
\end{proof}

\begin{corollary}\label{cor:qreg}
If $X$ is a quasiregular space with $F(X) \le \lambda$ and $(2^\lambda,\lambda)$ is a \pd\ for $X$
then $d(X) \le \lambda$.
\end{corollary}

If $X$ is any countably tight Lindelöf space then, as is well known, $F(X) \le \omega$ holds. Hence by Corollary
\ref{cor:qreg}  $X$ is separable if $(\mathfrak{c},\omega)$ is a \pd\ for $X$,
provided that $X$ is also quasiregular. Considering that
in \cite{BT} all spaces were assumed to be Tychonov, a condition much stronger than quasiregularity,
the above yields an affirmative answer to Problem 4.7, and hence also of 4.6 and 4.5 of \cite{BT},
the latter two being just particular cases of the first one.

Concerning 4.5 that is about countably tight $\sigma$-compact spaces we conjecture that such a ($T_3$-)space
$X$ already having $(\omega_1,\omega)$ -- instead of $(\mathfrak{c},\omega)$ -- as a \pd\ for $X$
implies that $X$ is separable.

We shall see later that the assumption of quasiregularity in Corollary \ref{cor:qreg} is essential because
it is consistent to have a non-separable $T_2$-space $X$ with $F(X) = \omega$ such that
$(\mathfrak{c},\omega)$ is even a \dom\ for $X$.

It was shown in 3.1 of \cite{GTW} that if $\big((2^\lambda)^+,\lambda\big)$ is a \dom\ for
a $T_3$-space $X$ then $d(X) \le \lambda$, and in Theorem 3.18 of \cite{BT} the same conclusion is
proved for Tychonov $X$ under the weaker assumption that
$\big((2^\lambda)^+,\lambda\big)$ is a \pd\ for $X$. In the proof of the latter the assumption of $X$ being
Tychonov was essential because a compactification of $X$ was used.
Our next result yields a significant improvement of both of these results.
First, however, we need a definition.

\begin{definition}\label{df:cross}
Let $X$ be a topological space and $\kappa$ ba a cardinal. We say that the family
$\{\langle U_\alpha,V_\alpha \rangle : \alpha < \kappa\}$ of {\em disjoint} pairs of open
subsets of $X$ is $\kappa$-{\em crossing} if for every set $I \in [\kappa]^\kappa$ of indices
there are $\alpha, \beta \in I$ such that $U_\alpha \cap V_\beta \ne \emptyset$.
\end{definition}

\smallskip

\begin{theorem}\label{tm:cross}
If $\kappa > 2^\lambda$ and $(\kappa,\lambda)$ is a \pd\ for the space $X$ then
$X$ does not admit a $\kappa$-crossing family.
\end{theorem}

\begin{proof}
Assume, on the contrary, that
$\{\langle U_\alpha,V_\alpha \rangle : \alpha < \kappa\}$ is a $\kappa$-crossing family of $X$ and
fix $S \in [X]^\lambda$ that pins down the $\le \kappa$-sized collection $\{U_\alpha \cap V_\beta : \alpha, \beta < \kappa\} \setm \{\emptyset\}$.
Since $S$ only has $2^\lambda < \kappa$ subsets, then there are $S_0, S_1 \subs S$
and $I \in [\kappa]^\kappa$ such that for all $\alpha \in I$ we have $S \cap U_\alpha = S_0$ and $S \cap V_\alpha = S_1$.
But then $U_\alpha \cap V_\alpha = \emptyset$ implies $S_0 \cap S_1 = \emptyset$ contradicting that for some
$\alpha, \beta \in I$ we have $S \cap U_\alpha \cap V_\beta \ne \emptyset$.
\end{proof}

At this point we recall that, for any space $X$, $RO(X)$ forms a complete boolean algebra and as such it admits an independent subset
$\mathcal{I}$ of cardinality $\varrho(X) = |RO(X)|$ by the celebrated result of Balcar and Franek in \cite{BF}.
But if $\mathcal{I} = \{U_\alpha : \alpha < \varrho(X)\}$ is an independent subset of $RO(X)$ and $V_\alpha = X \setm \overline{U_\alpha}$
then clearly $\{\langle U_\alpha,V_\alpha \rangle : \alpha < \varrho(X)\}$ is a $\varrho(X)$-crossing family of $X$, immediately
implying the following result.

\begin{corollary}\label{cor:ro}
If $\big((2^\lambda)^+,\lambda\big)$ is a \pd\ for any space $X$ then $\varrho(X) = |RO(X)| \le 2^\lambda$.
Consequently, we have $d(X) \le \lambda$ if $RO^+(X)$ is a $\pi$-base of $X$.
\end{corollary}

Note that obviously $RO^+(X)$ is a $\pi$-base of $X$ if $X$ is either semiregular or qusairegular.
This shows that Corollary \ref{cor:ro} really is much stronger than the above mentioned
analogous results from \cite{GTW} resp. \cite{BT}.

Let us recall next that for every $T_2$-space $X$ we have $|X| \le 2^{\varrho(X)}$ because every point
of $X$ is the intersection of its regular open neighborhoods. Consequently, if we also have that
$\big((2^\lambda)^+,\lambda\big)$ is a \pd\ for $X$ then Corollary \ref{cor:ro} implies $|X| \le 2^{2^\lambda}$.
This leads us to the following result.

\begin{corollary}\label{cor:t2}
If $\big(2^{2^\lambda},\lambda\big)$ is a \pd\ for a $T_2$-space $X$ then $d(X) \le \lambda$.
\end{corollary}

\begin{proof}
Indeed, Corollary \ref{cor:ro} implies $\varrho(X) \le 2^\lambda$, hence we have $|X| \le 2^{2^\lambda}$.
But then Theorem \ref{tm:card} can be applied with $\kappa = 2^{2^\lambda}$ to conclude $d(X) \le \lambda$.
\end{proof}

For $\lambda = \omega$ we thus have that if $(2^\mathfrak{c}, \omega)$ is a \pd\ for a $T_2$-space $X$ then
$X$ is separable. This yields an affirmative answer to question 4.2 of \cite{GTW}, which asked this
under the stronger condition that $(2^\mathfrak{c}, \omega)$ is a \dom\ for $X$. If $2^\mathfrak{c} = \mathfrak{c}^+$
then this also yields a consistent affirmative answer to question 4.1 of \cite{GTW}.

By Corollary \ref{cor:t2} the case $2^{2^\lambda} = (2^\lambda)^+$ of GCH is a sufficient condition for the validity of the statement
"if $\big((2^\lambda)^+,\lambda\big)$ is a \dom\ for a $T_2$-space $X$ then $d(X) \le \lambda$".
Next we show that this condition is also necessary.
This, of course, yields a consistent negative answer to question 4.1 of \cite{GTW}.

\begin{example}\label{ex:++}
If $2^{2^\lambda} > (2^\lambda)^+$ then there is an Urysohn space $X$ such that $\big((2^\lambda)^+,\lambda\big)$ is a \dom\ for
$X$ but $d(X) > \lambda$.
\end{example}

\begin{proof}
Let us put $\mu = 2^{2^\lambda}$ and $\nu = (2^\lambda)^{++}$, then $\mu = |\mathbb{C}_{2^\lambda}| \ge \nu$.
It is well known that the Cantor cube $\mathbb{C}_{2^\lambda}$ has a dense subset of size $\lambda$ but we claim that
actually it has $\mu$ many pairwise disjoint such dense subsets, say $\{D_\xi : \xi < \mu\}$.
The simplest way to see this is perhaps by noting that $\mathbb{C}_{2^\lambda}$ is actually a topological group
and so the dense subset $D$ of size $\lambda$ may be assumed to be a subgroup. Then the cosets by $D$
form a family of $\mu$ pairwise disjoint dense sets of size $\lambda$ in $\mathbb{C}_{\mu}$.

Let us now put $X = \bigcup \{D_\xi : \xi < \nu\}$ and fix a well-ordering $\prec$ of $X$ in order type $\nu$
such that for any $\xi < \zeta < \nu$ we have $D_\xi \prec D_\zeta$.
Then we consider the topology
$\tau$ on $X$ generated by its subspace topology inherited from $\mathbb{C}_{2^\lambda}$ and all final segments of $\prec$. The space $X$ obtained
in this way is Urysohn  because $\tau$ is finer than the subspace topology.

Now, if $Y$ is any proper initial segment of $\prec$ then,  on one hand, $Y$ is $\tau$-closed,
and on the other hand there is a $\xi < \nu$ such that $Y \prec D_\xi$ and hence
$Y$ is included in the $\tau$-closure of $D_\xi$. But this clearly implies both
$d(\langle X,\tau \rangle) = \nu > \lambda$ and that $\big((2^\lambda)^+,\lambda\big)$ is a \dom\ for
$\langle X,\tau \rangle$.
\end{proof}

Trivially, if $(\kappa,\lambda)$ is a \pd\ for a space $X$ then every regular cardinal
$\varrho$ satisfying $\lambda < \varrho \le \kappa$ is a caliber of $X$.
The "generic left separated" space $X$ on $\kappa = cf(\kappa) > \omega_1$ constructed in \cite{JSz}
is a hereditarily Lindelöf 0-dimensional space for which all regular cardinals
strictly between $\omega$ and $\kappa = \mathfrak{c}$ are calibers of $X$ but
$(\mathfrak{c},\omega)$ is not a \pd\ for $X$.

The authors of \cite{BT}
mentioned that they were unaware of the existence of a ZFC example of a space $X$ that has
$\omega_1$ as a caliber but $(\mathfrak{c},\omega)$ is not a \pd\ for $X$.
In what follows, we shall produce such ZFC examples.

\begin{definition}\label{df:tall}
Let $\kappa$ be a cardinal.
\begin{enumerate}[(i)]
\item
We say that a proper ideal $\mathcal{I}$ on $\kappa$ is $\kappa$-{\em tall} if
for every set $S \in [\kappa]^\kappa$ there is $I \in \mathcal{I}$ such that
$I \subs S$ and $|I| = |S| =\kappa$.
\item For every family $\mathcal{S}$ of subsets of $\kappa$ we put
$$\Sigma_\mathcal{S}(\kappa) = \{x \in \mathbb{C}_\kappa : supp(x) \in \mathcal{S}\},$$
where, of course, $supp(x) = \{\alpha < \kappa : x(\alpha) = 1\}$ is the support of $x$.
\end{enumerate}
\end{definition}

If $\kappa$ is an uncountable regular cardinal then the non-stationary ideal $NS(\kappa)$
is a $\kappa$-complete and $\kappa$-tall ideal with $[\kappa]^{< \kappa} \subs NS(\kappa)$.
Another such ideal on $\kappa$ is the one $\kappa$-generated by a maximal almost disjoint
family $\mathcal{A} \subs [\kappa]^{< \kappa}$ with $|\mathcal{A}|> \kappa$.
These examples of ideals show that the assumption of our next result is not vacuous.

\begin{theorem}\label{tm:tall}
Let $\kappa$ be an uncountable regular cardinal and $\mathcal{I}$ be
a $\kappa$-complete and $\kappa$-tall ideal with $[\kappa]^{< \kappa} \subs \mathcal{I}$.
Then $\kappa$ is a caliber of $\Sigma_\mathcal{I}(\kappa)$ but
$(\kappa,\lambda)$ is not a \pd\ for $\Sigma_\mathcal{I}(\kappa)$ for any $\lambda < \kappa$.
\end{theorem}

\begin{proof}
Let us consider a family $\{[\varepsilon_i] : i < \kappa\}$ of elementary open sets in $\mathbb{C}_{\kappa}$,
where $\{\varepsilon_i : i < \kappa\} \subs Fn(\kappa,2)$. We may assume without loss of generality that the
$\varepsilon_i$'s have the same size and that they form a $\Delta$-system
with root $\varepsilon$. Let $\eta_i = \varepsilon_i \setm \varepsilon$,
then $D_i = \{\alpha : \eta_i(\alpha) = 1\}$ may be written as
$D_i = \{\alpha_{i,k} : k < n\}$ for some $n < \omega$ and the $D_i$'s are pairwise disjoint.

Now, applying that $\mathcal{I}$ is $\kappa$-tall there is some $I_0  \in \mathcal{I}$ such
that $I_0  \subs \{\alpha_{i,0} : i < \kappa\}$ and $|I_0| = \kappa$. Let us put $K_0 = \{i < \kappa : \alpha_{i,0} \in I_0\}$,
then $|K_0| = \kappa$  implies that we can similarly find
$I_1  \in \mathcal{I}$ with  $|I_1| = \kappa$ such
that $I_1  \subs \{\alpha_{i,1} : i \in K_0\}$. Continuing this in $n$ steps we obtain sets
$J \in \mathcal{I} \cap [\kappa]^\kappa$ and $K \in [\kappa]^\kappa$  such that $\{\alpha_{i,k} : i \in K\} \supset  J$ for all $k < n$.
But then the for point $x \in \mathbb{C}_{\kappa}$ with $supp(x) = J \cup D$, where $D = \{\alpha : \varepsilon(\alpha) = 1\}$, we have
that $$x \in \bigcap \{[\varepsilon_i] : i \in K\} \ne \emptyset,$$
proving that $\kappa$ is indeed a caliber of $\Sigma_\mathcal{I}(\kappa)$.

On the other hand, the family $\{[\varepsilon_i] : i < \kappa\}$ of elementary open sets where
$supp(\varepsilon_i) = \{i\}$ clearly cannot be pinned down by $\lambda < \kappa$ points of
$\Sigma_\mathcal{I}(\kappa)$ because $\mathcal{I}$ is proper and $\kappa$-complete,
i.e. $(\kappa,\lambda)$ is not a \pd\ for $\Sigma_\mathcal{I}(\kappa)$.
\end{proof}

Applying Theorem \ref{tm:tall} for $\kappa = \omega_1$ and $\mathcal{I} = NS(\omega_1)$ we obtain the following
result that we promised above.

\begin{corollary}\label{cor:w1}
The subspace  $\Sigma_{NS}(\omega_1)$ of the Cantor cube $\mathbb{C}_{\omega_1}$ is a ZFC example of a space that has
$\omega_1$  as a caliber but for which $(\omega_1,\omega)$, and so $(\mathfrak{c},\omega)$ as well, is not a \pd .
\end{corollary}

In Proposition 3.11 of \cite{GTW} it was noted that if $X$ is $\lambda$-monolitic, i.e. $nw(\overline{Y}) \le \lambda$
for all $Y \in [X]^\lambda$, moreover $(2^\lambda,\lambda)$ is a \dom\ for $X$ then $d(X) \le \lambda$.
(In fact, having $(\lambda^+,\lambda)$  a \dom\ for $X$ suffices for this.) We'll show next that in this
result \dom\ may not be weakened to \pd .

\begin{example}\label{ex:mon}
For every cardinal $\lambda$ there is a Urysohn space $X$ in which every subset of size $\lambda$
is closed, hence $X$ is $\lambda$-monolithic, $(\kappa,\lambda)$ is a \pd\ for $X$
whenever $\kappa < 2^{2^\lambda}$, but $d(X) = \lambda^+$.
\end{example}

\begin{proof}
Our space $X$ is $\mathbb{C}_{2^\lambda}$ with the finer, hence Urysohn topology
$$\tau = \{U \setm A : U \in \tau(\mathbb{C}_{2^\lambda}) \text{ and } A \in [\mathbb{C}_{2^\lambda}]^{\le \lambda}\}.$$
We have seen in the proof of \ref{ex:++} that there are $2^{2^\lambda}$ pairwise disjoint $\lambda$-sized dense sets
$\{D_\xi : \xi < 2^{2^\lambda}\}$  in $\mathbb{C}_{2^\lambda}$.
So, if $\kappa < 2^{2^\lambda}$ and $\mathcal{U} = \{U_\alpha \setm A_\alpha : \alpha < \kappa\} \subs \tau^+(X)$ then
$|\bigcup \{A_\alpha : \alpha < \kappa\}| \le \kappa \cdot \lambda < 2^{2^\lambda}$
implies that there is
a $\xi < 2^{2^\lambda}$ such that $D_\xi \cap \bigcup \{A_\alpha : \alpha < \kappa\} = \emptyset$.
But we have $|U_\alpha| = 2^{2^\lambda}$ for all $\alpha < \kappa$,
hence $D_\xi$ clearly pins down $\mathcal{U}$.
Thus $(\kappa, \lambda)$ is indeed a \pd\ for $X$, while $d(X) = \lambda^+$ is obvious.
\end{proof}

We do not know if the Urysohn property in Example \ref{ex:mon} can be improved to $T_3$ or
Tychonov, so we consider that it yields only a partial answer to Problem 4.11 of \cite{BT} that is
asking if an $\omega$-monolitic space for which $(\mathfrak{c},\omega)$ is a \pd\ is
necessarily separable.

\section{Results for \dom s}

The aim of this section is to present solutions to a couple
of problems concerning \dom s that were posed in \cite{GTW}.
The first one yields a negative answer to Problem 4.3 of \cite{GTW}.

\begin{example}\label{ex:bQ}
There is a non-separable Urysohn space $X$ of countable $\pi$-character such
that $(\mathfrak{c},\omega)$ is a \dom\ for $X$.
\end{example}

\begin{proof}
The underlying set of $X$ is $\beta Q$, the \v Cech-Stone compactification of the rationals.
The topology $\tau$ of $X$ is the one generated by $\tau(\beta Q)$ together with all final
segments of a fixed well-ordering $\prec$ of the set $\beta Q$ of order type $2^\mathfrak{c}$.

Then $X$ is Urysohn since $\tau$ is finer than $\tau(\beta Q)$. Clearly, $X$ is non-separable
because $d(X) \ge cf(\prec) > \mathfrak{c}$.

Since $\beta Q$ is regular, we clearly have $\pi(\beta Q) = \pi(Q) = \omega$. This implies
that for every point $x \in X$ we have $\pi\chi(x,X) = \pi\chi(x,\beta Q) = \omega$ because
if $\{U_n : n < \omega\}$ is a local $\pi$-base of $x$ in $\beta Q$ and $V_n = \{y \in U_n : x \prec y\}$
then clearly $\{V_n : n < \omega\}$ is a local $\pi$-base of $x$ in $X$. So, $X$ has countable $\pi$-character.

To see that $(\mathfrak{c},\omega)$ is a \dom\ for $X$, let us note that if $x \in X$, $S \subs X$
with $x \prec S$ and $x$ is in the $\beta Q$-closure of $S$ then $x$ is also in the $X$-closure of $S$.
But it is known that $\Delta(\beta Q) = |\beta Q| = 2^\mathfrak{c}$ and this
together with $\pi(\beta Q) = \omega$ implies that every final segment of $\prec$ has a countable
subset that is dense in $\beta Q$. Putting these facts together we may then conclude that
$(\mathfrak{c},\omega)$ is indeed a \dom\ for $X$.
\end{proof}

Problem 4.10 of \cite{GTW} asked if spaces of countable pseudocharacter and having
$(\mathfrak{c},\omega)$ as a \dom\ have cardinality at most $2^\mathfrak{c}$. Our
next result says that actually there are 0-dimensional such spaces of arbitrarily
large cardinality.

\begin{example}\label{ex:psi}
For every cardinal $\kappa$ satisfying $\kappa^\mathfrak{c} = \kappa$ there is a dense
subspace $X$ of $\mathbb{C}_\kappa$ of countable pseudocharacter such that
$(\mathfrak{c},\omega)$ is a \dom\ for $X$ .
\end{example}

Before giving the proof of the statement of Example \ref{ex:psi} we present a lemma that
might be of independent interest.

\begin{lemma}\label{lm:0-1}
There is a countable dense subset $D$ of $\mathbb{C}_\mathfrak{c}$ such that
for every $\alpha < \mathfrak{c}$ and $x \in D$ we have
$$|\{n < \omega : x(\alpha + n) = 1\}| < \omega.$$
\end{lemma}

\begin{proof}
For every ordinal $\alpha < \mathfrak{c}$ let us consider the half-open
interval $$B_\alpha = \big[\omega \cdot \alpha,(\omega \cdot \alpha) + \omega \big),$$ i.e. the $\alpha$th
"block" of length $\omega$ in $\mathfrak{c}$ and take the product space $$P = \prod \{[B_\alpha ]^{< \omega} : \alpha < \mathfrak{c}\},$$
where each factor $[B_\alpha]^{< \omega}$ carries the (countable) discrete topology.

It is known, see e.g. 5.4 of \cite{Ju}, that $P$ has a countable dense subset, say $S$.
For every point $s \in S$ we then define a point $x_s \in \mathbb{C}_\mathfrak{c}$ with the
following stipulation: $$x_s(\omega \cdot \alpha + n) = 1\,\Leftrightarrow\,n \in s(\alpha)$$
and then put $D = \{x_s : s \in S\}$. Now, it is straight forward to check that $D$ is as required.
\end{proof}

\begin{proof}[Proof of \ref{ex:psi}]
Similarly as in the proof of Lemma \ref{lm:0-1}, for all $\alpha < \kappa$ we shall use $B_\alpha$ to denote the $\alpha$th
block $\big[\omega \cdot \alpha,(\omega \cdot \alpha) + \omega \big)$ of length $\omega$, but now
in $\kappa$, not in $\mathfrak{c}$. For any subset $S$ of $\kappa$ we also put $B_S = \bigcup \{B_\alpha : \alpha \in S\}$.

Using that $\kappa^\mathfrak{c} = \kappa$,  and, in particular that $cf(\kappa) > \mathfrak{c}$,
we may define by an easy transfinite recursion of length $\kappa$ a function $F : [\kappa]^\mathfrak{c} \to \kappa$
such that
\begin{enumerate}[(i)]
\item $B_S < \{F(S)\}$ for every $S \in [\kappa]^\mathfrak{c}$;

\smallskip

\item if $S,T \in [\kappa]^\mathfrak{c}$ and $S \ne T$ then $F(S) \ne F(T)$.
\end{enumerate}
We then put $R = \ran (F) = \{F(S) : S \in [\kappa]^\mathfrak{c}\}$ and will define the subspace
$X$ of $\mathbb{C}_\kappa$ that we are looking for in the form $X = \{x_\xi : \xi \in B_R\}$.

To do that, for any $S \in [\kappa]^\mathfrak{c}$ we may apply Lemma \ref{lm:0-1} to the copy
$\{0,1\}^{B_S}$  of $\mathbb{C}_\mathfrak{c}$ and define points $x_{\xi(n)} \in \mathbb{C}_\kappa$ with indices
in the block $$B_{F(S)} = \{\xi(n) = F(S) + n : n < \omega\}$$
by the following stipulations:
\begin{enumerate}[(a)]
\item $\{x_{\xi(n)}\upharpoonright B_S : n < \omega\}$ is dense in $\{0,1\}^{B_S}$ such that, for any $\alpha \in S$ and $n < \omega$,
we have $x_{\xi(n)}(\zeta) = 1$ only for finitely many $\zeta \in B_\alpha$;

\smallskip

\item for any $n,k < \omega$ we have $x_{\xi(n)}(\xi(k)) = 0$ iff $k < n$;

\smallskip

\item for every $n < \omega$ we have $x_{\xi(n)}(\zeta) = 0$ for all $\zeta \in \kappa \setm B_S \cup B_{F(S)}$.
\end{enumerate}
Since $F(S) \ne F(T)$ implies $B_{F(S)} \cap B_{F(T)} = \emptyset$, this can be done without any conflict.
Note that for every index $\xi \in B_R$ the point $x_\xi$ satisfies $|supp(x_\xi)| \le \mathfrak{c}$,
hence we actually have $$X = \{x_\xi : \xi \in B_R\} \subs \Sigma_\mathfrak{c}(\kappa).$$

Now, it is obvious that we have $|X| = |B_R| = \kappa$ and that $X$ is dense in $\mathbb{C}_\kappa$.
Also, for every set $A \in [X]^\mathfrak{c}$ there is $S \in [\kappa]^\mathfrak{c}\}$ such that
for every point $x \in A$ we have $supp(x) \subs B_S$. But then it easily follows from our construction that $A$ is contained in the closure
of the countable set $\{x_\xi : \xi \in B_{F(S)}\}$, hence $(\mathfrak{c},\omega)$ is a \dom\ for $X$.

To see that $X$ has countable pseudocharacter, fix  $\alpha \in R$ and $n < \omega$ and then consider the elementary $G_\delta$-set
$H_{\alpha,n}$ in $\mathbb{C}_\kappa$ with support $B_\alpha$ defined by
$$H_{\alpha,n} = \{x \in \mathbb{C}_\kappa : x(\alpha + k) = 0\, \Leftrightarrow\, k < n\}.$$
Now, using that $supp(x_\xi) \cap B_\alpha$ is finite for every $\xi \in B_R \setm B_\alpha$, it is easy to check that
$X \cap H_{\alpha,n} = \{x_{\alpha + n}\}$, hence $\psi(x_{\alpha + n},X) = \omega$,
which completes the proof.
\end{proof}

\section{Two forcing constructions}

In this section we shall present two consistency results via forcing concerning \dom s.
We put them in a separate section not just because they only yield consistency results
but also because they are much more technical than the rest of the paper.

Our first construction shows that Corollary \ref{cor:qreg} is not provable for $T_2$-spaces, even if
the condition of $(2^\lambda,\lambda)$ being a \pd\ is strengthened to being a \dom\
for a $T_2$-space.

\begin{theorem}\label{tm:T2}
It is consistent with CH that there is a non-separable $T_2$-space $X$ such that $(\mathfrak{c},\omega) = (\omega_1,\omega)$
is a \dom\ for $X$ and $F(X) = \omega$.
\end{theorem}

\begin{proof}
The basic idea of our construction is to obtain by forcing a $T_2$ topology $\tau$ on $\omega_2 = \mathfrak{c}^+$
such that, denoting by $cl_\tau(C)$ the $\tau$-closure of a set $C$, we have
\begin{enumerate}[(i)]
\item
for every $\alpha < \omega_2$ the block $B_\alpha = \big[\omega \cdot \alpha,(\omega \cdot \alpha) + \omega \big)$ is $\tau$-dense, and

\smallskip

\item $cl_\tau(C) \cap cl_\tau(D) \cap \omega_1 \ne \emptyset$ for any two infinite subsets $C$ and $D$ of $\omega_2$.
\end{enumerate}

Once we have got $\tau$, we pass to the finer topology $\rho$ on $\omega_2$ that is generated by $\tau$ together with
the final segments $[\alpha,\omega_2)$ of $\omega_2$ but only for $\alpha \ge \omega_1$.
Still, it is obvious that then $X = \< \omega_2,\rho \>$ has density $\omega_2$.
It is clear by (i) that for any bounded subset $S$ of $\omega_2$ with $S < B_\alpha$ we have
$S \subs cl_\rho(B_\alpha)$, hence $(\mathfrak{c},\omega)$ is a \dom\ for $X$.
Finally, it is also clear that (ii) remains valid for $\rho$, consequently every free sequence in $X$
has order-type $< \omega \cdot 2$, hence we have $F(X) = \omega$.

Now, it "only" remains to obtain $\tau$ as above. To do that, we
start with a ground model $V$ that satisfies CH, and we are going to define a notion of forcing $\mathbb{P}$
that will produce this topology $\tau$.

$\mathbb{P}$  will be both countably closed and will satisfy the $\omega_2$-CC, hence cardinals and cofinalities will be  preserved in
the generic extension $V^\mathbb{P}$. Moreover, CH will hold in $V^\mathbb{P}$ and no new countable subsets
of $V$ will be added.

To define $\mathbb{P} = \<P, \le \>$, we first fix an arbitrary index set $I$ of cardinality $\omega_1$.
We shall call a subset $A \subs \omega_2$ {\em saturated} if $A \cap B_\alpha \ne \emptyset$ implies $B_\alpha \subs A$
for any $\alpha < \omega_2$. In other words, the saturated sets are just unions of blocks.

The elements of $P$ will be of the form $p = \<U_p,\mathcal{C}_p\>$, where $U_p$ is a function with
domain $I_p \in [I]^\omega$ and values taken in $[\omega_2]^\omega$ in such a way that $$A_p = \bigcup \{U_p(i) : i \in I_p\}$$
is saturated.
There is more that we demand of $U_p$ however. Namely, we have requirements about the topology $\tau_p$ on $A_p$
that is generated by the family $\{U_p(i) : i \in I_p\}$ as a subbase. (Of course, $\tau_p$ is an approximation
of our desired topology $\tau$.)

For any $J \in [I_p]^{< \omega}$ we let $V_{p,J} = \bigcap \{U_p(i) : i \in J\}$. Clearly, then
$$\mathcal{V}_p = \big\{V_{p,J} \ne \emptyset : J \in [I_p]^{< \omega}\big\}$$ forms a base for the topology $\tau_p$.
What we require then is that for each block $B_\alpha \subs A_p$ and
$V_{p,J} \in \mathcal{V}_p$ we have $$|B_\alpha \cap V_{p,J}| = \omega.$$
In particular, then $B_\alpha$ is $\tau_p$-dense.
It should be mentioned here that we do not require $\tau_p$ to satisfy any separation axiom.

The second coordinate $\mathcal{C}_p$ of the condition $p$ is a countable collection of pairs
$\langle x,C \rangle \in A_p \times [A_p]^\omega$ such that $x \in acc_p(C)$, i.e. $x$ is a $\tau_p$-accumulation point of $C$.
Equivalently, this means that $x \in V \in \mathcal{V}_p$ implies $|V \cap C| = \omega$.

Our next step is to define for $p,q \in P$ when $p \leq q$, i.e. $p$ is a stronger condition than $q$.

\begin{definition}\label{df:leq1}
For $p,q \in P$ we define $p \leq q$ by the following stipulations:
\begin{enumerate}[(a)]
\item $I_p \supset I_q$;

\smallskip

\item $U_q(i) = A_q \cap U_p(i)$ for every $i \in I_q$;

\smallskip

\item $V_{q,J} \cap V_{q,K} = \emptyset$ implies $V_{p,J} \cap V_{p,K} = \emptyset$ for any two $J, K \in [I_q]^{< \omega}$;

\smallskip

\item $\mathcal{C}_p \supset \mathcal{C}_q$.
\end{enumerate}
\end{definition}
We note that (a) and (b) trivially imply $A_p \supset A_q$. Also, we emphasize that if $p \leq q$ then
for every pair $\langle x,C \rangle \in \mathcal{C}_q$ we must have $x \in acc_p(C)$, that is not
automatic from $x \in acc_q(C)$.

We next present several lemmas concerning the forcing $\mathbb{P}$ which together will yield us the desired
topology $\tau$ on $\omega_2$.

\begin{lemma}\label{lm:wcl1}
The forcing $\mathbb{P}$ is countably closed.
\end{lemma}

\begin{proof}[Proof of \ref{lm:wcl1}]
Assume that $p_0 \ge p_1 \ge ...$ is a decreasing $\omega$-sequence in $\mathbb{P}$, where
$p_n = \< U_n, \mathcal{C}_n \>$ for $n < \omega$. (So, to enhance readability, we use $n$ instead of $p_n$
as indices.) We shall define a lower bound $p = \<U_p,\mathcal{C}_p\>$ for the $p_n$'s as follows.

\begin{enumerate}
\item $I_p = \bigcup_{n < \omega} I_n$;

\smallskip

\item $U_p(i) = \bigcup_{n \ge n_i} U_n(i)$ for $i \in I_p$, where  $n_i = \min \{n : i \in I_n\}$;

\smallskip

\item $\mathcal{C}_p = \bigcup_{n < \omega} \mathcal{C}_n$.
\end{enumerate}

We first check that $p \in P$. For this, note that (1) and (2) imply $A_p = \bigcup_{n < \omega} A_n$,
moreover $V_{p,J} = \bigcup_{n \ge n_J} V_{n,J}$ for $J \in [I_p]^{< \omega}$, where  $n_J = \min \{n : J \subs I_n\}$.
It immediately follows that for the topology $\tau_p$ generated by $\{U_p(i) : i \in I_p\}$ on $A_p$ we have
$\tau_n = \tau_p \upharpoonright A_n$ for all $n < \omega$, moreover if
$V_{p,J} \in \mathcal{V}_p$ then we have $|B_\alpha \cap V_{p,J}| = \omega$ for any $B_\alpha \subs A_p$.

We still have to show that $x \in acc_p(C)$ for any $\langle x,C \rangle \in \mathcal{C}$. Using the above and (3), if
$x \in V_{p,J}$ then there is $n < \omega$ for which we have both $x \in V_{n,J}$ and $\langle x,C \rangle \in \mathcal{C}_n$,
so $|C \cap V_{n,J}| = \omega$ that implies $|C \cap V_{p,J}| = \omega$. Thus we indeed have $p \in P$.

To check that $p \le p_n$ for each $n < \omega$, only item (c) needs a little thought. But we have $p_m \le p_n$  for any $m > n$,
hence $V_{n,J} \cap V_{n,K} = \emptyset$ implies $V_{m,J} \cap V_{m,K} = \emptyset$ for all $m > n$,
consequently we have $V_{p,J} \cap V_{p,K} = \emptyset$ as well.
\end{proof}

Our next lemma is an amalgamation result that will be used in the proof of the fact that $\mathbb{P}$ satisfies the $\omega_2$-CC.
First we need a definition. We shall say two conditions $p,q \in P$ are {\em isomorphic} provided that the following hold.

\begin{enumerate}
\item $I_p = I_q$;

\smallskip

\item $A = A_p \cap A_q < A_p \setm A_q  < A_q \setm  A_p$;

\smallskip

\item $otp(A_p) = otp(A_q)$ and if $\varphi$ is the unique order isomorphism taking $A_p$ onto $A_q$
(which is the identity on $A$) then we have $U_q(i) = \varphi[U_p(i)]$  for all $i \in I_p = I_q$;

\smallskip

\item $\mathcal{C}_q = \{\langle \varphi(x), \varphi[C]\rangle : \< x,C \> \in \mathcal{C}_p\}$.
\end{enumerate}

\begin{lemma}\label{lm:CC1}
If $p,q \in P$ are isomorphic then they are compatible in $\mathbb{P}$.
\end{lemma}

\begin{proof}[Proof of \ref{lm:CC1}]
Let us define the condition $r = \< U_r, \mathcal{C}_r \> \in P$ by setting $U_r(i) = U_p(i) \cup U_q(i)$
for $i \in I_p$, moreover $\mathcal{C}_r = \mathcal{C}_p \cup \mathcal{C}_q$. Clearly, we have then
that $A_r = A_p \cup A_q$ is saturated.

Of course, we have to check that $r \in P$. It is clear from the definition of $r$ that for any finite $J \subs I_r = I_p$
we have $V_{r,J} = V_{p,J} \cup V_{q,J}$. Hence $\tau_p = \tau_r \upharpoonright A_p$ and $\tau_q = \tau_r \upharpoonright A_q$,
moreover $\tau_r \upharpoonright A = \tau_p \upharpoonright A = \tau_q \upharpoonright A$. It is clear then that for every
$V_{r,J} \in \mathcal{V}_r$ and $B_\alpha \subs A_r$ we have $|B_\alpha \cap V_{p,J}| = \omega.$
Also, if $\langle x,C \rangle \in \mathcal{C}_p$ then by the above $x \in acc_p(C)$ implies $x \in acc_r(C)$ and the same is
true with $q$ replacing $p$.

Next, we have to show that$r$ extends both $p$ and $q$. For this only item (c) of \ref{lm:CC1} might not look immediately obvious.
So, assume e.g. that $V_{p,J} \cap V_{p,K} = \emptyset$ for $J, K \in [I_p]^{< \omega}$. Then by the isomorphism of $p$ and $q$
we have $V_{q,J} \cap V_{q,K} = \emptyset$, hence $V_{r,J} \cap V_{r,K} = \emptyset$  as well. This shows $r \le p$ and
$r \le q$ is shown analogously.
\end{proof}

Now, using CH, for any $S \in [P]^{\omega_2}$ a standard counting and $\Delta$-system argument yields two members of $S$
that are isomorphic and hence compatible, demonstrating that $\mathbb{P}$ satisfies the $\omega_2$-CC.

Our next lemmas yield various extension properties of $\mathbb{P}$ which will be used in establishing the desired
properties of the generic topology $\tau$.

\begin{lemma}\label{lm:Aq}
For every condition $q \in P$ and $\alpha \in \omega_2$ with $B_\alpha \cap A_q = \emptyset$ there is $p \le q$
such that $B_\alpha \subs A_p$.
\end{lemma}

\begin{proof}[Proof of \ref{lm:Aq}]
Let us fix an enumeration $A_q = \{a_k : k < \omega\}$ of $A_q$ and a partition $\{\Gamma_k : k < \omega\}$
of $B_\alpha$ into infinite sets.
We then define $p = \<U_p,\mathcal{C}_p\>$ such that $I_p = I_q$ by letting $$U_p(i) = U_q(i) \cup \bigcup \{\Gamma_k : a_k \in  U_q(i)\} $$
for all $i \in I_q$, and
$\mathcal{C}_p = \mathcal{C}_q$. Then $\tau_q = \tau_p \upharpoonright A_q$ is obvious, as well as $|V \cap B_\alpha| = \omega$
for all $V \in \mathcal{V}_p$. Also, then $x$ is a $\tau_p$-accumulation point of $C$ whenever
$\langle x,C \rangle \in \mathcal{C}_p = \mathcal{C}_q$. Thus indeed $p \in P$.

To check condition (c) for $p \le q$ observe that for every $V_{q,J} \in \mathcal{V}_q$ we have
$$V_{p,J} = V_{q,J} \cup \bigcup \{\Gamma_k : a_k \in  V_{q,J}\},$$ hence
$V_{q,J} \cap V_{q,K} = \emptyset$ implies $V_{p,J} \cap V_{p,K} = \emptyset$.
\end{proof}

\begin{lemma}\label{lm:Iq}
For every condition $q \in P$ and $j \in I \setm I_q$ there is $p \le q$
such that $j \in I_p$.
\end{lemma}

\begin{proof}[Proof of \ref{lm:Iq}]
It is straight forward to check that if we put $U_p = U_q \cup \{\<j,A_q\>\}$ then
$p = \<U_p,\mathcal{C}_q\>$ is as required.
\end{proof}

Our following lemma will be used to conclude that the generic topology $\tau$ is $T_2$.

\begin{lemma}\label{lm:T2q}
For every condition $q \in P$ and $\{x_0,x_1\} \in [A_q]^2$ there is $p \le q$
such that there are disjoint $E_0,E_1 \in \mathcal{V}_p$ with $x_0 \in E_0$ and $x_1 \in E_1$.
\end{lemma}

\begin{proof}[Proof of \ref{lm:T2q}]
Let us consider the following countable family of infinite subsets of $A_q$:
$$\mathcal{A}_q = \{V \cap B_\alpha : V \in \mathcal{V}_q \text{ and } B_\alpha \subs A_q\} \cup \{V \cap C : x \in V \in \mathcal{V}_q \text{ and }\<x,C\> \in \mathcal{C}_q \}.$$
Then $\mathcal{A}_q$ can be reaped, i.e. there is a partition $A_q = E_0 \cup E_1$ of $A_q$ such that both $E_0$ and $E_1$ have infinite
intersection with all members of $\mathcal{A}_q$. Of course, we may also assume that $x_0 \in E_0$ and $x_1 \in E_1$.

Let us fix two new indices $j_0, j_1 \in I \setm I_q$,
then define $U_p$ with domain $I_p = I_q \cup \{j_0,j_1\}$ by $U_p(i) = U_q(i)$ for $i \in I_q$,
moreover $U_p(j_0) = E_0$ and $U_p(j_1) = E_1$. Now, put $p = \<U_p,\mathcal{C}_q\>$,
then it is easy to check that both $p \in P$ and $p \le q$ are satisfied, completing the proof.
\end{proof}

The next lemma is needed to establish property (ii) of the generic topology $\tau$.

\begin{lemma}\label{lm:w1q}
Given any  condition $q \in P$ and two infinite subsets $C,D \subs A_q$, there are $p \le q$ and
$x \in \omega_1 \cap A_p$ such that both $\<x,C\> \in \mathcal{C}_p$ and $\<x,D\> \in \mathcal{C}_p$.
\end{lemma}

\begin{proof}[Proof of \ref{lm:w1q}]
Let us pick an index $j \in I \setm I_q$ and an ordinal $\alpha < \omega_1$ such that $B_\alpha \cap A_q = \emptyset$.
Also, fix a partition $\{\Gamma_k : k < \omega\}$ of $B_\alpha$ into infinite sets with $x = \omega \cdot \alpha \in \Gamma_0$.
We also fix an enumeration $A_q = \{a_k : k < \omega\}$ of $A_q$.

Then we define $U_p$ with domain $I_p = I_q \cup \{j\}$ as follows: $U_p(j) = A_q \cup \Gamma_0$ and
$$U_p(i) = U_q(i) \cup \bigcup \{\Gamma_{k+1} : a_k \in  U_q(i)\}$$ whenever $i \in I_q$. Finally,
we put $\mathcal{C}_p = \mathcal{C}_q \cup \{\<x,C\> ,\<x,D\> \}$.

To see that $p = \<U_p,\mathcal{C}_p\> \in P$ note first that $\tau_q$ is the subspace topology of $\tau_p$
on $A_q$, hence for any $\<y,E\> \in \mathcal{C}_q$ we have $y \in acc_p(E)$. It is also clear that
$x \in V_{p,J}$ implies $J = \{j\}$, hence $V_{p,J} = U_p(J) \supset A_q \supset C \cup D$.
Consequently, we have both $x \in acc_p(C)$ and $x \in acc_p(D)$.

Finally, to check condition (c) for $p \le q$, we observe that $$V_{p,J} = V_{q,J} \cup \bigcup \{\Gamma_{k+1} : a_k \in  V_{q,J}\}$$
for any $J \in [A_q]^{<\omega}$. It immediately follows then that
$V_{q,J} \cap V_{q,K} = \emptyset$ implies $V_{p,J} \cap V_{p,K} = \emptyset$ for any two $J, K \in [I_q]^{< \omega}$.
\end{proof}

Having these lemmas, we are now ready to obtain the desired topology $\tau$ in the generic extension $V^\mathbb{P}$.
First, however we note that Lemmas \ref{lm:wcl1} and \ref{lm:CC1} imply that $\omega_1$ and $\omega_2$ are preserved and CH holds in $V^\mathbb{P}$.

So, let $G$ be $\mathbb{P}$-generic over our ground model $V$ and we define the map $U : I \to \mathcal{P}(\omega_2)$
in the generic extension $V[G]$ of $V$ by the formula
$$U(i) = \bigcup \{U_p(i) : p \in G\}.$$ It follows from Lemma \ref{lm:Iq} that $U$ is indeed defined for all $i \in I$, while
Lemma \ref{lm:Aq} on the other hand implies that $\bigcup \{U(i) : i \in I\} = \omega_2$.

Now, $\tau$ is the topology on $\omega_2$ generated by the range $\{U(i) : i \in I\}$ of the map $U$ as a subbase. So, putting
$V_J = \bigcap \{U(i) : i \in J\}$ we obtain a base $$\mathcal{V} = \{V_J : J \in [\omega_2]^{< \omega}\} \setm \{\emptyset\}$$ for $\tau$.
Note that for any such $J$ we have $$V_J = \bigcup \{V_{p,J} : p \in G \text{ and } J \subs I_p\}.$$

It is immediate from the definition of the conditions $p \in P$ and Lemma \ref{lm:Aq} that every block $B_\alpha$ is $\tau$-dense, which
was requirement (i) for $\tau$. On the other hand, Lemma \ref{lm:T2q} implies that $\tau$ is $T_2$. Indeed, if $x_0,x_1$ are distinct
points in $\omega_2$ then the set $D_{x_0,x_1}$ of those $p \in P$ for which
there are disjoint $E_0,E_1 \in \mathcal{V}_p$ with $x_0 \in E_0$ and $x_1 \in E_1$ is dense in  $\mathbb{P}$ by Lemma \ref{lm:T2q}.
But clearly, if $p \in G \cap D_{x_0,x_1}$ then $p$ forces that $x_0$ and $x_1$ have disjoint $\tau$-neighborhoods.
This is where requirement (c) in Definition \ref{df:leq1} is needed.

Finally, we show that requirement (ii) is also satisfied by $\tau$. This is actually immediate from Lemma \ref{lm:w1q} by which,
for any two countably infinite subsets $C$ and $D$ of $\omega_2$, the set of $p \in P$ with some $x \in \omega_1 \cap A_p$ and
both $\{\<x,C\>, \<x,D\>\} \subs \mathcal{C}_p$ is dense in $\mathbb{P}$. So, there is such a $p \in G$, while $p$ clearly
forces that $$x \in \omega_1 \cap acc_\tau(C) \cap acc_\tau(D) \ne \emptyset.$$
This completes the proof of the Theorem.
\end{proof}

We would like to point out that the above dense-in-itself $T_2$-space $X$ of cardinality $> \mathfrak{c}$ has the property that
any two infinite closed subsets of it intersect, i.e. it is strongly anti-Urysohn using the terminology from \cite{JSSz2}.
So, as a collateral corollary we obtained by this a consistent answer to a problem of \cite{JSSz2} that asked if
strongly anti-Urysohn spaces of cardinality $> \mathfrak{c}$ could exist. However, as it happens, in a very recent
and still unwritten joint work with S. Shelah
we could actually produce such a space in ZFC.

\medskip

     Problem 4.6 of \cite{GTW} raised the following question: If $(\mathfrak{c},\omega)$
is a \dom\ for a Fr\`echet-Urysohn (in short: FU) space $X$, is then $X$ separable? Our next forcing construction yields a consistent
negative answer to this question.

\begin{theorem}\label{tm:FU}
It is consistent with CH that there is a non-separable 0-dimensional $T_2$-space $X$ such that $(\mathfrak{c},\omega) = (\omega_1,\omega)$
is a \dom\ for $X$ and $X$ is FU.
\end{theorem}

\begin{proof}
The space $X$ that we shall construct will be left-separated in order type $\omega_2$, i.e. all its initial segments in the left-separating
well-ordering will be closed. This clearly will imply $d(X) = \omega_2$.
Before starting our forcing construction, we need some terminology that, naturally, concerns left-separation.

Let $S$ be any set of ordinals, considered with its natural well-ordering. We say that the map $U : S \to \mathcal{P}(S)$ is left-separating
(in short: LS) if for every $\xi \in S$ we have $\min U(\xi) = \xi$. Any such LS map $U$ determines a 0-dimensional $T_2$
topology $\tau_U$ on $S$ that has $$\{U(\xi) : \xi \in S\} \cup \{S \setm U(\xi) : \xi \in S\}$$ as a subbase and is left-separated
by its natural well-ordering.

Given such an LS map $U : S \to \mathcal{P}(S)$, for any $\xi \in S$ we put $U^1(\xi) = U(\xi)$ and $U^0(\xi) = S \setm  U(\xi)$, and
for every finite function $\varepsilon \in Fn(S,2)$ with domain $D(\varepsilon)$
we  let $$[\varepsilon]_U  = \bigcap \{U^{\varepsilon(\xi)}(\xi) : \xi \in D(\varepsilon)\}.$$
Clearly, then $\mathcal{V}_U = \{[\varepsilon]_U : \varepsilon \in Fn(S,2)\}$ is a base of $\tau_U$.
Our aim then is to apply
forcing to obtain a model of CH in which there is an LS map $U$ on some $X \subs \omega_2$ with $|X| = \omega_2$ such that
(a) $(\omega_1,\omega)$ is a \dom\ for $\tau_U$ and (b) $\tau_U$ is FU.

To achieve this aim,
we start with a ground model $V$ that satisfies CH, and then define a countably closed and $\omega_2$-CC notion of forcing $\mathbb{P}$
that will produce the required LS map $U$. Cardinals and cofinalities will be  preserved in
the generic extension $V^\mathbb{P}$, moreover, CH will hold in $V^\mathbb{P}$ and no new countable subsets
of $V$ will be added.
Having the goals (a)  and (b) in mind, the somewhat complicated definition of our notion of forcing $\mathbb{P} = \< P, \le \>$
should seem quite natural.

Also, as several parts of our new forcing construction is somewhat parallel with the earlier one,
in what follows we shall freely use the terminology and notation that was introduced there.
A new piece of notation that we shall need is the following. Given a saturated set of ordinals $A$
and a block $B_\alpha \subs A$, we write $A/\alpha$ to denote $A \cap (\omega \cdot \alpha + \omega)$.
In other words, $A/\alpha$ is the initial segment of $A$ that has $B_\alpha$ as its top block.

\begin{definition}\label{df:P2}
First we define the set of conditions $P$. The elements of $P$ will be pairs of the form $p = \< U_p, \mathcal{C}_p \>$,
where $U_p$ is an LS map on a countable saturated set $A_p \subs \omega_2$ and  $\mathcal{C}_p$ is a
countable subset of $A_p \times [A_p]^\omega$, so that the following two requirements are satisfied.
(As in the previous construction, we write $\tau_p$ instead of $\tau_{U_p}$, etc.)

\begin{enumerate}[(i)]
\item If $B_\alpha \subs A_p$ and $V \cap A_p/\alpha \ne \emptyset$ for some $V \in \tau_p$ then
$|V \cap B_\alpha| = \omega$.

\smallskip

\item If $\langle x,C \rangle \in \mathcal{C}_p$ then $C$ converges to $x$ in $\tau_p$ (in symbols: $C \to_p x$), i.e.
we have $C \subs^* V\,(\Leftrightarrow\, |C \setm V| < \omega)$ whenever $x \in V \in \tau_p$.
\end{enumerate}

\smallskip

Next we define for $p,q \in P$ when $p \leq q$, i.e. $p$ is a stronger condition than $q$.
So, $p \leq q$ iff

\begin{enumerate}[(1)]
\item $A_p \supset A_q$,

\smallskip

\item $U_p(\xi) \cap A_q = U_q(\xi)$ for all $\xi \in A_q$,

\smallskip

\item $\mathcal{C}_p \supset \mathcal{C}_q$.
\end{enumerate}
It is easy to see that $\mathbb{P} = \< P, \le \>$ is indeed a partial ordering.
\end{definition}

Let us make a few comments about this definition. First, (i) means that any block $B_\alpha \subs A_p$
is on one hand $\tau_p$-dense in $A_p/\alpha$, moreover $B_\alpha$ is dense-in-itself in $\tau_p$.
Also, when checking (i), it suffices to consider the countably many basic $\tau_p$-open sets $V \in \mathcal{V}_p$.
On the other hand, to check if $C \to_p x$ for some $x \in A_p$ and $C \in [A_p]^\omega$, it suffices to prove that $x \in {U_p}(\xi)$
implies $C \subs^* {U_p}(\xi)$ and $x \in A_p \setm  {U_p}(\xi)$
implies $|C \cap {U_p}(\xi)| < \omega$. This observation comes in handy
when, in establishing the relation $p \leq q$, we have to check that for any $\langle x,C \rangle \in \mathcal{C}_q$
we have not just $C \to_q x$ but even $C \to_p x$.

We next will formulate and prove several lemmas that will establish the required properties of $\mathbb{P}$.
The first lemma says that $\mathbb{P}$ is countably closed in a strong sense.

\begin{lemma}\label{lm:wcl2}
Any decreasing $\omega$-sequence in $\mathbb{P}$ has a greatest lower bound,
hence $\mathbb{P}$ is countably closed.
\end{lemma}

\begin{proof}[Proof of \ref{lm:wcl2}]
Assume that $p_0 \ge p_1 \ge ...$ is a decreasing $\omega$-sequence in $\mathbb{P}$, where
$p_n = \< U_n, \mathcal{C}_n \>$ for $n < \omega$. We then define $p = \< U_p, \mathcal{C}_p \>$
by putting $A_p = \bigcup_{n < \omega} A_n$ and $U_p(\xi) = \bigcup_{n_\xi \le n < \omega} U_n(\xi)$
for all $\xi \in A_p$, where  $n_\xi = \min \{n : \xi \in A_n\}$,
moreover $\mathcal{C}_p = \bigcup_{n < \omega} \mathcal{C}_n$.

To prove that $p \in P$ and $p \le p_n$ for each $n < \omega$ is very similar to the proof of Lemma
\ref{lm:wcl1} and so we leave it to the reader to check the details. It is  also obvious that $p$
is the greatest lower bound of the $p_n$'s.
\end{proof}

Now we start to work toward showing that $\mathbb{P}$ is $\omega_2$-CC. Just like in the previous construction,
this will follow from the fact that isomorphic conditions in $\mathbb{P}$ can be amalgamated but here this
will require much more work. We actually start with presenting a purely combinatorial result.

\begin{lemma}\label{lm:id}
Let $\mathcal{I}$ be a proper ideal on $\omega$ with $[\omega]^{< \omega} \subs \mathcal{I}$ and assume that we
are given three countable families:
$\mathcal{V} \subs \mathcal{I}^+ = \mathcal{P}(\omega) \setm \mathcal{I},$
$\mathcal{D} = \{D_n : n < \omega\} \subs \mathcal{I}$, and $\mathcal{C} = \{\< x_n, C_n \> : n < \omega\} \subs \omega \times \mathcal{I}$
such that

\begin{enumerate}[(a)]
\item $|C_n \cap D_m| < \omega$ for any $n,m \in \omega$;

\smallskip

\item $x_n \ne x_m$ implies $|C_n \cap C_m| < \omega$.
\end{enumerate}
Then there is a set $E \subs \omega$ for which
\begin{enumerate}[(i)]
\item $E$ reaps $\mathcal{V}$, i.e. $|V \cap E| = |V \setm E| = \omega$ for all $V \in \mathcal{V}$,

\smallskip

\item $x_n \in E$ implies $C_n \subs^* E$ and $x_n \notin E$ implies $|C_n \cap E| < \omega$,

\smallskip

\item  $|D_n \cap E| < \omega$ for all $n < \omega$.
\end{enumerate}
\end{lemma}

\begin{proof}[Proof of \ref{lm:id}]
We start by taking an $\omega$-abundant enumeration $\mathcal{V} = \{V_n : n < \omega\}$, i.e. such that
for every $V \in \mathcal{V}$ we have $|\{n < \omega : V_n = V\}| = \omega$. Then
we are going to define by recursion sets $E_n, F_n \in \mathcal{I}$ such that $E_n \cap F_n = \emptyset$,
moreover, for each $n < \omega$, the following three inductive assumptions are satisfied:

\begin{enumerate}[(1)]
\item $\{x_k : k < n\} \subs E_n \cup F_n$,

\smallskip

\item $E_n =^* \bigcup \{C_k : k < n \text{ and } x_k \in E_n\}$,

\smallskip

\item $F_n =^* \bigcup \{C_k : k < n \text{ and } x_k \in F_n\} \cup \bigcup \{D_k : k < n \}$.
\end{enumerate}
(Of course, $A = ^* B$ means that the symmetric difference $A \triangle B$ of $A$ and $B$ is finite.)

For $n = 0$ we may simply put $E_0 = F_0 = \emptyset$. Assuming that $E_n, F_n \in \mathcal{I}$ with $E_n \cap F_n = \emptyset$
have been defined and the inductive hypotheses hold, we shall define $E_{n+1}$ and $F_{n+1}$.
First, as $W_n = V_n \setm (E_n \cup F_n \cup \{x_n\} \cup C_n \cup D_n) \in \mathcal{I}^+$ is infinite, we pick two distinct
elements $y_n,z_n \in W_n$.

Now we distinguish two cases. If $x_n \in F_n$ then we observe that both $C_n \cap E_n$ and $D_n \cap E_n$ are finite by
our assumptions and so we may put $E_{n+1} = E_n \cup \{y_n\}$ and $F_{n+1} = F_n \cup \{z_n\} \cup C_n \cup D_n \setm E_n$.
If, on the other hand $x_n \notin F_n$, then we use that both $C_n \cap F_n$ and $D_n \cap (E_n \cup C_n)$ are finite, hence
we may put $E_{n+1} = E_n \cup \{y_n\} \cup \{x_n\} \cup C_n\setm F_n$ and $F_{n+1} = F_n \cup \{z_n\} \cup D_n \setm E_{n+1}$.
It is simple to check that the inductive assumptions are satisfied for $n+1$ in both cases.

Having completed this recursive construction, we claim that $E = \bigcup \{E_n : n < \omega\}$ satisfies the requirements
of our lemma. Indeed, the choice of the $y_n$'s and $z_n$'s clearly implies (i), i.e. that $E$ reaps $\mathcal{V}$.
Moreover, the validity of our inductive assumptions immediately implies (ii) and (iii).
\end{proof}

We next apply Lemma \ref{lm:id} to obtain a topological result that will be crucial in proving our desired amalgamation result.
But first we need a definition.

\begin{definition}\label{df:sui}
We call a topology $\tau$ on a countable saturated set $A$ of ordinals {\em suitable} if it is 0-dimensional, $T_2$, second countable,
moreover $B_\alpha \subs A\,$ and $V \in \tau$ with $V \cap A/\alpha \ne \emptyset$ implies $|V \cap B_\alpha| = \omega$.
\end{definition}
Clearly, $\tau_p$ is suitable whenever $p \in P$ but the converse of this fails. Still, this more general concept will play an important role below.

\begin{lemma}\label{lm:sui}
Let $\tau$ be a suitable topology on a countable saturated set $A$ of ordinals left-separated
by its natural well-ordering
and $T$ be a (possibly empty) saturated proper initial segment of $A$.
Let us also be given a countable collection $\mathcal{K} \subs A \times [A]^\omega$ such that $C \to_\tau x$
for every $\langle x,C \rangle \in \mathcal{K}$.

Then there is a subset $E$ of the final segment $S = A \setm T$ of $A$ for which
the topology $\varrho$ on $A$ generated by $\tau \cup \{E,A \setm E\}$ is suitable,
both $B_\alpha \cap E$ and $B_\alpha \setm E$ are $\tau$-dense in $B_\alpha$  whenever $B_\alpha \subs S$,
moreover $C \to_\varrho x$ for all $\langle x,C \rangle \in \mathcal{K}$.
\end{lemma}

\begin{proof}[Proof of \ref{lm:sui}]
Since $|S| = \omega$, we may "transfer" Lemma \ref{lm:id} from $\omega$ to $S$ by considering the ideal $\mathcal{I}$  on $S$
whose members are those subsets of $S$ that can be written as the union of finitely many $\tau$-discrete sets.
Clearly, for every $\tau$-converging set $C$ we have $C \cap S \in \mathcal{I}$, moreover $B_\alpha \subs S\,$ and $V \in \tau$
with $V \cap A/\alpha \ne \emptyset$ imply $V \cap B_\alpha \in \mathcal{I}^+$ because $V \cap B_\alpha$ is $\tau$-dense-in-itself,
while all members of $\mathcal{I}$ are $\tau$-scattered.

Let $\mathcal{B}$ be a countable base of $\tau$ and then we define the following three countable families:
\begin{align*}
      \mathcal{V} &= \{V \cap B_\alpha \ne \emptyset : V \in \mathcal{B} \text{ and } B_\alpha \subs S\} \subs \mathcal{I}^+,\\
      \mathcal{D} &= \{C \cap S : \langle x,C \rangle \in \mathcal{K} \text{ and } x \in T\},\\
      \mathcal{C} &= \{\< x,C \cap S\> : \langle x,C \rangle \in \mathcal{K} \text{ and } x \in S\}.
\end{align*}
Now, it is easy to check that they satisfy all conditions of Lemma \ref{lm:id} "transferred" to $S$.
Note that the left separation of $\tau$ implies that $C \subs^*  S$ if $\langle x,C \rangle \in \mathcal{K} \text{ and } x \in S.$

We claim that the set $E \subs S$ given by Lemma \ref{lm:id} will be as required. The topology $\varrho$ is clearly suitable,
moreover the fact that $E$ reaps $\mathcal{V}$ means that $B_\alpha \cap E$ and $B_\alpha \setm E$ are
$\tau$-dense in $B_\alpha$  for all $B_\alpha \subs S$.

If $\langle x,C \rangle \in \mathcal{K}$ with $x \in T \subs A \setm E$ then we have $C \subs^* A \setm E$.
If, on the other hand, $\langle x,C \rangle \in \mathcal{K}$ with $x \in S$ then $x \in E$ implies $C \subs^* E$,
while $x \in S \setm E$ implies $C \subs^* S \setm E \subs  A \setm E$. Clearly, then we have
$C \to_\varrho x$ for all $\langle x,C \rangle \in \mathcal{K}$.
\end{proof}

Let us point out that the new topology $\varrho$ on $A$ together with $\mathcal{K}$ satisfies all the assumptions of Lemma
\ref{lm:sui}, hence it can be applied to it again. This then allows us to iterate this process through all the natural numbers
and thus yields the following immediate corollary of Lemma \ref{lm:sui}.

\begin{corollary}\label{cor:wsui}
Let $A, T, S, \tau$ and $\mathcal{K}$ be as in Lemma \ref{lm:sui}. Then we may find subsets $\{E_n : n < \omega\}$ of $S$
such that for each $n < \omega$ the topology $\tau_n$ generated by $\tau \cup \{E_k : k < n\} \cup \{A \setm E_k : k < n\}$
is suitable, $C \to_{\tau_n} x$ for all $\langle x,C \rangle \in \mathcal{K}$, moreover for each $n > 0$ and $B_\alpha \subs S$
both $B_\alpha \cap E_n$ and $B_\alpha \setm E_n$ are $\tau_{n - 1}$-dense in $B_\alpha$.
\end{corollary}

We are now ready to present the amalgamation result that will be used to prove that $\mathbb{P}$ is $\omega_2$-CC.
We call two conditions $p,q \in P$ are {\em isomorphic} provided that the following hold.

\begin{enumerate}
\item $A_p \cap A_q < A_p \setm A_q  < A_q \setm  A_p$;

\smallskip

\item $otp(A_p) = otp(A_q)$ and if $\varphi$ is the unique order isomorphism from $A_p$ to $A_q$
(which is the identity on $A_p \cap A_q$) then we have $U_q(i) = \varphi[U_p(i)]$  for all $i \in I_p = I_q$;

\smallskip

\item $\mathcal{C}_q = \{\langle \varphi(x), \varphi[C]\rangle : \< x,C \> \in \mathcal{C}_p\}$.
\end{enumerate}

\begin{lemma}\label{lm:CC2}
If $p,q \in P$ are isomorphic then they are compatible in $\mathbb{P}$. Consequently, $\mathbb{P}$ is $\omega_2$-CC.
\end{lemma}

\begin{proof}[Proof of \ref{lm:CC2}]
The common extension $r = \< U_r,\mathcal{C}_r \>$ of $p$ and $q$ will be chosen so that the domain of $U_r$ is $A_r = A_p \cup A_q$,
moreover $\mathcal{C}_r = \mathcal{C}_p \cup \mathcal{C}_q$.
Then we are actually forced to define the values of $U_r$ for $x \in A_p \cap A_q$ by $U_r(x) = U_p(x) \cup U_q(x)$
and for $x \in A_q \setm A_p$ by $U_r(x) = U_q(x)$.

The non-trivial part is to define $U_r(x)$ for $x \in A_p \setm A_q$. To do that, we fix a one-one enumeration
$A_p \setm A_q = \{a_n : n < \omega\}$ and apply Corollary \ref{cor:wsui} for the suitable topology $\tau_q$ on $A_q$,
its initial segment $T = A_p \cap A_q$ and the family $\mathcal{K} = \mathcal{C}_q$.
The latter gives us subsets $\{E_n : n < \omega\}$ of $A_q \setm  A_p$ and then we put $U_r(a_n) = U_p(a_n) \cup E_n$
for each $n < \omega$.

First of all, we have to check that $r \in P$. Condition (i) only needs to be verified for $B_\alpha \subs A_q \setm  A_p$. This, however
follows from the fact that all the topologies $\tau_n$ from Corollary \ref{cor:wsui} are suitable. This same same Corollary also
yields that if $\langle x,C \rangle \in \mathcal{C}_q$ then $C \to_{\tau_r} x$, while $\tau_p = \tau_r \upharpoonright A_p$
trivially implies that if $\langle x,C \rangle \in \mathcal{C}_p$ then $C \to_{\tau_r} x$, so we verified condition (ii) as well.
But from this and the definition of $r$ it is obvious that $r \le p,q$.

Finally, using CH and standard counting and $\Delta$-system arguments we can show that among any $\omega_2$ elements of $P$
there are two that are isomorphic, hence compatible. Consequently, $\mathbb{P}$ is indeed $\omega_2$-CC.
\end{proof}

Let us now consider any $\mathbb{P}$-generic filter $G$ over our ground model $V$ and, in $V[G]$, the LS map $U$
with domain $X = \bigcup \{A_p : p \in G\}$ defined by the stipulation $$U(x) = \bigcup \{U_p(x) : p \in G \text { and } x \in A_p\}.$$
Clearly, $X$ has a canonical $\mathbb{P}$-name $\dot{X}$ and for any condition $p \in P$ if $x \in A_p$
then $p \Vdash x \in \dot{X}$. Also, it is obvious that $X \subs \omega_2$ is saturated. Our next result implies that
$X$ is also cofinal in $\omega_2$, hence $|X| = \omega_2$.

\begin{lemma}\label{lm:cof}
For every condition $p \in P$ and ordinal $\xi < \omega_2$ there is $r \le p$ such that $A_r \setm \xi\ne \emptyset$.
So, $\mathbb{P} \Vdash |X| = \omega_2$.
\end{lemma}

\begin{proof}[Proof of \ref{lm:cof}]
We may assume that $A_p \subs \xi$ and consider an isomorphic copy $q$ of $p$ such that $\xi \le \min A_q$.
Then we may use Lemma \ref{lm:CC2} to obtain a common extension $r$ of  $p$ and $q$ and then $r$ is as required.
\end{proof}

While it is obvious that $x \in A_p$ implies $p \Vdash x \in \dot{X}$, the converse of this is not.
What we do have is the following result.

\begin{lemma}\label{lm:xin}
If $p \Vdash x \in \dot{X}$ then there is an extension $r \le p$ such that $x \in A_r$.
Equivalently, $p \Vdash x \in \dot{X}$ iff $\{r \le p : x \in A_r\}$ is dense below $p$.
\end{lemma}

\begin{proof}[Proof of \ref{lm:xin}]
Let $G$ be any $\mathbb{P}$-generic filter with $p \in G$. Then $p \Vdash x \in \dot{X}$ implies $x \in X$,
hence there is $q \in G$ with $x \in A_q$ by the definition of $X$. But $G$ is a filter, hence there is
$r \in G$ that is a common extension of $p$ and $q$, and $r$ is as required.
\end{proof}

We are now ready to return to the proof of our actual goal, namely that the LS topology $\tau_U$ on $X$ satisfies both
(a) $(\omega_1,\omega)$ is a \dom\ for $\tau_U$ and (b) $\tau_U$ is FU.

Now, (a) clearly follows if we can show that each $B_\alpha \subs X$ is dense in $X/\alpha$. To see this, we have to show that
$[\varepsilon]_U \cap X/\alpha \ne \emptyset$ for any basic open set $[\varepsilon]_U \in \mathcal{V}_U$ implies
$[\varepsilon]_U \cap B_\alpha \ne \emptyset$. So, pick $x \in [\varepsilon]_U \cap X/\alpha$ and use Lemma \ref{lm:xin} and the
countable closedness of $\mathbb{P}$ to obtain $p \in G$ such that $D(\varepsilon) \cup B_\alpha \cup \{x\} \subs A_p$.
But condition (i) for $p \in P$ then yields $[\varepsilon]_p \cap B_\alpha \ne \emptyset$, hence by $[\varepsilon]_p \subs [\varepsilon]_U$
we are done.

Next we start working towards our goal (b). First we observe that if $p \in G$ and  $\<x,C\> \in \mathcal{C}_p$ then conditions
(ii) and (3) from Definition \ref{df:P2} trivially imply that $C \to_{\tau_U} x$. In other words, this means that
$p \Vdash C \to_{\tau_U} x$ for any $p \in P$ and $\<x,C\> \in \mathcal{C}_p$. Consequently, it suffices to prove that whenever
$\dot{Z}$ is a $\mathbb{P}$-name for a subset of $X$ and $p \Vdash x \in acc_U(\dot{Z})$ then $p$ has an extension $r \le p$
such that for some $C \in [A_r]^\omega$ we have $r \Vdash C \subs \dot{Z}$ and $\<x,C\> \in \mathcal{C}_r$.

To see this, assume that
$p \Vdash x \in acc_U(\dot{Z})$. We may assume $x \in A_p$ by Lemma \ref{lm:xin}. Then for any $\varepsilon \in Fn(A_p,2)$ with
$x \in [\varepsilon]_p$ we have $p \Vdash \dot{Z} \cap [\varepsilon]_U \setm \{x\} \ne \emptyset$. Now, a straightforward closure argument, repeatedly using Lemma \ref{lm:xin} and the
countable closedness of $\mathbb{P}$, yields us $q \le p$ such that for every $\varepsilon \in Fn(A_q,2)$ with
$x \in [\varepsilon]_q$ there is $z \in A_q \setm \{x\}$ such that $q \Vdash z \in \dot{Z} \cap [\varepsilon]_U$. Since $q \Vdash "\tau_q = \tau_U \upharpoonright A_q"$,
this implies that $x \in acc_{\tau_q}(Y)$, where $Y = \{z \in A_q : q \Vdash z \in \dot{Z}\}$.

But the topology $\tau_q$ is second countable and 0-dimensional, so it is metrizable and hence FU, consequently there is $C \in [Y]^\omega$
such that $C \to_{\tau_q} x$. But then $r =\<U_q,\mathcal{C}_q \cup \{\<x,C\>\}\>$ is an extension of $q$ and hence of $p$
that clearly is as required.
\end{proof}

We conclude by formulating the following natural questions that are left open.

\begin{problem}
\begin{enumerate}
\item Is there a ZFC example of a non-separable ($T_2$ or $T_3$) FU space $X$ such that $(\mathfrak{c},\omega)$
is a \dom\ for $X$?

\smallskip

\item Is it consistent to have a ($T_2$ or $T_3$) FU space $X$ of cardinality $> 2^\mathfrak{c}$ such that $(\mathfrak{c},\omega)$
is a \dom\ for $X$?
\end{enumerate}
\end{problem}

\end{document}